\documentclass[leqno]{amsart}
\usepackage[margin=1in]{geometry}
\usepackage{mathtools}
\usepackage{enumerate}
\usepackage{amsthm}
\usepackage{mathabx}
\usepackage[all]{xy}
\usepackage[titletoc,toc,title]{appendix}
\usepackage[outdir=./]{epstopdf}
\usepackage{graphicx,epstopdf} 
\usepackage{comment}
\usepackage{tikz-cd}
\newtheorem{theorem}{Theorem}[section] 

\newtheorem{lem}[theorem]{Lemma}

\newtheorem*{cl}{Claim}

\theoremstyle{definition}

\newtheorem*{rem}{Remark}

\theoremstyle{remark}
\newtheorem*{notation}{Notation}
\newtheorem*{ack}{Acknowledgment}
\newtheorem{example}{Example}[section]

\DeclareMathOperator{\interior}{\text{int}}
\DeclareMathOperator{\supp}{\textup{supp}}

\numberwithin{equation}{section}

\title{A topological space associated to corank 1 tropical phased matroids}

\author{Uly Alvarez}

\address{\noindent Ulysses Alvarez, Department of Mathematics,
The University of Alabama,
Tuscaloosa, AL 35487, USA}
\email{uaalvarez@ua.edu}

\date{}

\keywords{}

\begin{document}
\fontsize{12}{13pt} \selectfont  
    
\begin{abstract}
A consequence of the Topological Representation Theorem in \cite{FL} is that the geometric realization of the order complex
of the poset of non-zero covectors of a loopless rank $n-1$ oriented matroid on $[n]$ is homeomorphic to an $(n-2)$-sphere.
In this paper, we begin the study of an analogous theorem for tropical phased matroids by proving that the topological order complex for 
a loopless rank $n-1$ tropical phased matroid on $[n]$ is homeomorphic to a $(2n-3)$-sphere.
\end{abstract}

\maketitle

\section{Introduction}

Tropical phased matroids are examples of \textit{matroids over hyperfields} \cite{BB}, where a \textit{hyperfield} is a field with a multivalued addition.
A tropical phased matroid is a matroid over the \textit{tropical phase hyperfield} $\Phi$,
$\Phi$ is the union of the unit circle $S^1$ and the origin $\{0\}$ in the complex plane $\mathbb{C}$ 
equipped with the complex multiplication and addition $\boxplus$ such that
\begin{enumerate}
\item $x\boxplus0=\{x\}$,
\item $x\boxplus-x=S^1 \cup \{0\}$ whenever $x \neq 0$ and
\item $x \boxplus y$ is the smallest closed arc in $S^1$ joining $x$ and $y$ when $y\neq-x$.
\end{enumerate}

The notion of matroids over hyperfields is a generalization of oriented matroids, which can be viewed as 
matroids over the \textit{sign hyperfield} $\mathbb{S}$, where $\mathbb{S}$ is $\{-1,0,1\}$ 
equipped with the multiplication inherited from $\mathbb{R}$ and addition $\boxplus$ such that 
\begin{enumerate}
\item $x\boxplus0=\{x\}$,
\item $1\boxplus-1=\{-1,0,1\}$ whenever $x \neq 0$ and
\item $x \boxplus x=\{x\}$.
\end{enumerate}

We begin the study of carrying over the following result from the theory of oriented matroids to the realm of tropical phased matroids.

\begin{theorem}[The Topological Representation Theorem in \cite{FL}]
Let $\mathcal{M}$ be a loopless rank $r$ oriented matroid.
Then the poset of nonzero covectors of $\mathcal{M}$ is isomorphic to the face poset of the cell decomposition induced by an arrangement of pseudospheres in $S^{r-1}$, the standard unit sphere of dimension $r-1$.
\end{theorem}

As a consequence, we have that the geometric realization of the order complex of the non-zero covectors of $\mathcal{M}$ is homeomorphic to $S^{r-1}$.
In this paper, we will prove the following.
\begin{theorem}\label{main}
The topological order complex of the space of non-zero covectors of a given loopless rank $n-1$ tropical phased matroid on $[n]$
is homeomorphic to $S^{2n-3}$. 
\end{theorem}

The topological order complex, first defined in \cite{Zi}, was not the first topological space that was considered for this endeavor. 
In \cite{AD}, the topological space considered arose from the search for a topology on $\Phi$ such that the \textit{phase} map
\begin{center}
ph$\colon \mathbb{C} \to \Phi,\  z \mapsto
\begin{cases}
z/|z| &\text{ if } z \neq 0,\\
0 & \text{ otherwise}
\end{cases}
$
\end{center}
is continuous, namely, the open sets of $\Phi$ are the usual open sets of $S^1$ together with the set $S^1 \cup \{0\}$.
Since the covectors of a tropical phased matroid on $[n]$ lives in $\Phi^n$, this topology on $\Phi$ induces a topology on the set of covectors.
However this space is significantly obscure whereas, by viewing the set of covectors as a topological poset, its topological order complex is more transparent.
For example, the former is not generally Hausdorff but the latter is. 
The following theorem allows us to choose to study the latter case.

\begin{theorem}[\cite{AG}]
Let $\mathcal{V^*}$ be the set of covectors of a tropical phased matroid.
The topological order complex of $\mathcal{V^*}$ is weak homotopy equivalent to $\mathcal{V^*}$ endowed with the topology considered in \cite{AD}.
\end{theorem}

\begin{ack}
I am indebted to Laura Anderson and Ross Geoghegan for the time and patience they invested in me, 
and for allowing me to grow under their guidance.

\end{ack}

\section{Background}

\subsection{Brief note on the piecewise linear category}
The bulk of this paper is concerned with polyhedra and their joins, so piecewise linear (PL) topology is the natural category for us.
However, it is often more convenient to name smooth ``round'' balls and their various intersections, rather than explicitly describing polyhedral versions.
The rule is that all the constructions in the paper can easily be seen to be homeomorphic to corresponding constructions
among rectilinear polyhedra in such a way that standard theorems of PL topology can be used.
For example,

\begin{lem}[\cite{RS}, Lemma 1.10]\label{extension}
Let $B_1$ and $B_2$ be $d$-balls and $h\colon \partial B_1 \to \partial B_2$ be a homeomorphism.
Then $h$ extends to a homeomorphism $B_1 \to B_2$.
\end{lem}

\begin{lem}[\cite{RS}, Corollary 3.16]\label{base}
Let $B_1$ and $B_2$ be $d$-balls such that 
\begin{enumerate}
\item $B_1 \cap B_2$ is a $(d-1)$-ball and
\item $B_1 \cap B_2=\partial(B_1) \cap \partial(B_2)$, where $\partial(B_k)$ is the boundary of $B_k$.
\end{enumerate}
Then $B_1 \cup B_2$ is a $d$-ball.
\end{lem}

The following result is the main tool of this paper, a generalization of Lemma \ref{base}.
\begin{theorem}\label{gluing}
Let $B_1$, $B_2$, ..., $B_m$ be $d$-balls, where $m \leq d+1$. 
If for every $J \subseteq [m]$ such that $|J|>1$,
\begin{enumerate}
\item $\bigcap\limits_{k \in J}B_k$ is a $(d-|J|+1)$-ball, and
\item for every $r \in J$, $\bigcap\limits_{k \in J} B_k \subset \partial\left(\bigcap\limits_{k \in J-\{r\}} B_k\right)$,
\end{enumerate}
then $\bigcup\limits_{k \in [m]} B_k$ is a $d$-ball. 
\end{theorem}
\begin{proof}[Sketch of proof]
We will induct on $m$. Lemma \ref{base} is precisely the case when $m=2$. 

Assume true for $m-1$. Let $\{B_k:k \in [m]\}$ be a collection of $d$-balls which satisfy the hypotheses.
Thus $\bigcup\limits_{k=1}^{m-1} B_k$ is a $d$-ball.
By Lemma \ref{base}, it suffices to show that  
$\left(\bigcup\limits_{k=1}^{m-1} B_k\right) \cap B_m$ is a $(d-1)$-ball which follows from hypothesis (1), and 
\begin{center}
$\left(\bigcup\limits_{k=1}^{m-1}B_k\right) \cap B_m=\partial\left(\bigcup\limits_{k=1}^{m-1} B_k\right) \cap \partial(B_m)$
\end{center}
which follows from hypothesis (2).
Thus $\left(\bigcup\limits_{k=1}^{m-1} B_k\right) \cup B_m$ is a $d$-ball.
\end{proof}

\subsection{The set of covectors as a topological poset}
To formally define the set of covectors of a given matroid over $\Phi$, one must first define the corresponding set of $\Phi$-\textit{circuits}.
However, since the general framework of matroids over hyperfields is not used in this paper, the set of covectors will simply be given:
let $\mathbf{v}=(v_1,\dots,v_n) \in (S^1)^n$ and define $$\mathbf{v}^\perp\coloneqq\{\mathbf{x}\in \Phi^n:0 \in v_1x_1 \boxplus \cdots \boxplus v_nx_n\}$$
to be the set of \textit{covectors} of a loopless rank $n-1$ tropical phased matroid on $[n]$.
In order to obtain the desired topological space, we endow $\mathbf{v}^\perp$ with a topological poset structure.

A \textit{topological poset} (\cite{Zi}) is a Hausdorff space $\mathcal{P}$ equipped with a partial ordering $\leq$ such that $\{(x,y): x \leq y \text{ in } \mathcal{P}\}$ is closed in $\mathcal{P}^2$. 
\begin{example}
Every \textit{discrete poset}, that is, a poset given the discrete topology, is a topological poset. 
\end{example}
\begin{example}
Endow $\Phi$ with a partial ordering such that the only non-trivial relation is that 0 is less than every point in $S^1$ and 
we equip it with the subspace topology inherited from $\mathbb{C}$. Then $\Phi$ is a topological poset. 
\end{example}
\begin{example}
The topological poset structure on $\Phi$ described in the previous example induces a topological poset structure on $\Phi^n$ 
by defining the partial ordering component-wise and taking the product topology. 
Thus, we have a topological poset structure on $\mathbf{v}^\perp-\{\mathbf{0}\}$ by restricting the partial ordering considering the subspace topology.
\end{example}

Given a topological poset $\mathcal{P}$, we can associate another space called 
the \textit{topological order complex} $\Delta(\mathcal{P})$ of $\mathcal{P}$.
For an arbitrary topological poset $\mathcal{P}$, its topological order complex is
the geometric realization of the simplicial space $\bigsqcup\limits_{k=0}^\infty N_k(\mathcal{P})$, 
where 
$$N_k(\mathcal{P})\coloneqq\{(x_0,x_1,\ldots,x_k) \in \mathcal{P}^{k+1}: x_0 \leq x_1 \leq \cdots \leq x_k \text{ in } \mathcal{P}\}.$$

There is a more digestible description of $\Delta(\mathcal{P})$ provided that $\mathcal{P}$ is a \textit{mirrored poset}, that is,
$\mathcal{P}$ is a topological poset with an associated a finite poset $Q$ and a poset map $\mu: \mathcal{P} \to Q$ such that 
\begin{enumerate}
\item for every $k \in Q$, $\mu^{-1}(k)$ is a non-empty and closed subset of $\mathcal{P}$ and 
\item $x < y$ in $\mathcal{P}$ implies $\mu(x) < \mu(y)$ in $Q$.
\end{enumerate}

For such a $\mathcal{P}$, $\Delta(\mathcal{P})$ is homeomorphic to
$$\left\{\sum_{k \in Q} t_k \mathbf{x}_k \in \bigast_{k \in Q}\mu^{-1}(k): 
\{\mathbf{x}_k: t_k \neq 0\} \text{ is a chain in } \mathcal{P}\right\},$$
where $\bigast_{k \in Q}\mu^{-1}(k)$ is the topological join of the family of spaces $\{\mu^{-1}(k): k \in Q\}$,
$t_k \in [0,1]$ for every $k \in Q$, and $\sum_{k \in Q} t_k=1$.
Informally speaking, $\Delta(\mathcal{P})$ can be obtained by taking $\mathcal{P}$ and gluing $k$-simplices along the $k$-chains of $\mathcal{P}$.

\begin{rem}
If $\mathcal{P}$ is a discrete poset, $\Delta(\mathcal{P})$ is the geometric realization of the order complex of $\mathcal{P}$
in the classical sense.
\end{rem}

\begin{example}
$\Phi$ is a mirrored poset by defining $\mu:\Phi \to \{0 < 1\}$ which maps 0 to 0 and all of $S^1$ to 1.
By definition,
$\Delta(\Phi)=\{0\} \ast S^1$ which is homeomorphic to the (closed) unit ball $\mathbb{B} \subset \mathbb{C}$ via 
\begin{center}
$g\colon \Delta(\Phi) \to \mathbb{B},\ (1-s)\cdot 0 + s\cdot x \mapsto sx$,
\end{center}
where $s\cdot x$ is the formal multiplication in a topological join and
$sx$ is the scalar multiplication in $\mathbb{C}$.
\end{example}

\begin{example}
$\Phi^n$ is a mirrored poset by defining
$\mu\colon \Phi^n \to [0,n]_\mathbb{Z}$ by $\mu(\mathbf{x})=|\supp(\mathbf{x})|$,
where $\supp(\mathbf{x}) \coloneqq \{k \in [0,n]_\mathbb{Z}: x_k \neq 0\}$ 
and $|\supp(\mathbf{x})|$ is the cardinality of $\supp(\mathbf{x})$.
Since the partial ordering is defined component-wise, $(x_1,\dots, x_n) < (y_1,\dots,y_n)$ in $\Phi^n$
implies $x_j=y_j$ for every $j \in \supp(\mathbf{x})$.
So if $\sum t_k \mathbf{x}_k \in \Delta(\Phi^n)$, we know what $\{\mathbf{x}_k: t_k \neq 0\}$ must satisfy.
However, to see that $\Delta(\Phi^n)$ is homeomorphic to $\mathbb{B}^n$, we need the following result.

\begin{lem}\label{homeo}
$\Delta(\Phi^n)$ and $(\Delta\Phi)^n$ are homeomorphic.
\end{lem}
\begin{proof}
Given an element $z=\sum_{k=0}^n t_k \mathbf{x}_k \in \Delta(\Phi^n)$,
let $x_{k,j}$ be the $j^{th}$-coordinate of $\mathbf{x}_k$.
Define $N(z,j) \coloneqq\min\{k \in [0,n]_\mathbb{Z} : t_k \neq 0 \text{ and } x_{k,j} \neq 0\}$ and
$$\ell(z,j) \coloneqq
\begin{cases}
N(z,j) &\text{ if $N(z,j)$ is non-empty},\\
n+1 & \text{ otherwise.}
\end{cases}
$$

Define $\gamma \colon\Delta(\Phi^n) \to (\Delta(\Phi))^n$ by 
$$\gamma(z) \coloneqq \left(\left(\sum_{k=0}^{\ell(z,j)-1} t_k\right) \cdot 0
+\left(\sum_{k=\ell(j,z)}^{n} t_k\right) \cdot x_{k,\ell(j,z)}\right)_{j=1}^n.$$
If $\ell(z,j)=n+1$, we assert $\left(\sum_{k=\ell(j,z)}^{n} t_k\right) \cdot x_{k,\ell(j,z)}=0$ .

When $\Phi$ is treated as a set, $\gamma$ is bijective \cite[Theorem 3.2]{Wa}.

To show that $\gamma$ is continuous, consider the diagram
\begin{center}
\begin{tikzcd}
\bigsqcup\limits_{k=0}^\infty N_k(\Phi^n) \arrow[r] \arrow[d, "q_n"]
& \left(\bigsqcup\limits_{k=0}^\infty N_k(\Phi)\right)^n \arrow[d, "q_1^n" ] \\
\Delta(\Phi^n) \arrow[r, "\gamma"]
&  (\Delta\Phi)^n
\end{tikzcd}
\end{center}
where $q_1$ and $q_n$ are the appropriate quotient maps, and
the top map is defined by $$(\mathbf{x}_0,\dots,\mathbf{x}_k) \mapsto ((x_{j,1})_{j=0}^k,\dots,(x_{j,n})_{j=0}^k).$$
Since the diagram commutes, $\gamma$ is continuous, and, thus, a homeomorphism.
\end{proof}
\end{example}

Furthermore, for every $y \in S^1$, we have
$g_y \colon\Phi^n \to \Phi^n$ defined by $g_y(\mathbf{x})=(yx_1,\ldots,yx_n)$. 
Clearly, this map preserves chains in $\Phi^n$, i.e., 
if $\{\mathbf{x}_0 < \cdots < \mathbf{x}_k\}$ is a chain in $\Phi^n$, then so is $\{g_y(\mathbf{x}_0) < \cdots < g_y(\mathbf{x}_k)\}$.
Thus  the induced map 
$$\Delta(g_y)\colon\Delta(\Phi^n) \to \Delta(\Phi^n),\ \sum t_k \mathbf{x}_k
\mapsto \sum t_k (g_y(\mathbf{x}_k))$$
is a homeomorphism.

\begin{example}\label{important}
For $\mathbf{v}=(v_1,\dots,v_n) \in (S^1)^n$, $\mathbf{v}^\perp-\{\mathbf{0}\}$ is a mirrored poset by restricting the domain of $\mu$
to $\mathbf{v}^\perp-\{\mathbf{0}\}$ and the codomain to $[2,n]_{\mathbb{Z}}$.
Note that 0 is not in the codomain since we removed $\mathbf{0}$ from the domain and 
1 is not in the codomain since there does not exist an element $\mathbf{x}=(x_1,\dots,x_n)$ such that 
$|\supp(\mathbf{x})|=1$ and $0 \in v_1x_1 \boxplus \cdots \boxplus v_nx_n$.
Moreover, we can define the homeomorphism
\begin{center}
$\mathbf{1}_n^\perp \to \mathbf{v}^\perp,\ (x_1,\dots,x_n) \mapsto (v_1x_1,\dots,v_nx_n)$
\end{center}
where $\mathbf{1}_n=(1,\dots,1) \in (S^1)^n$.
Though this map does not preserve chains, 
the preimage of every $k$-chain in $\mathbf{v}^\perp$ will be a unique $(k+1)$-set in $\mathbf{1}_n^\perp$.
Thus $\Delta(\mathbf{v}^\perp)$ is homeomorphic to $\Delta(\mathbf{1}_n^\perp)$ and, hence,
it suffices to understand the case when $\mathbf{v}=\mathbf{1}_n$.
\end{example}

\section{Understanding $\Delta(\mathbf{1}_n^\perp-\{\mathbf{0}\})$}

\subsection{When $0 \in x_1 \boxplus \cdots \boxplus x_n$}

We have addressed in Example \ref{important} why we are not interested in the cases when $n=0$ nor when $n=1$.

\begin{lem}\label{zero} Let $\mathbf{x}=(x_{1},\dots, x_{n}) \in\Phi^{n}-\{\bf 0\}$, where $|\supp(\mathbf{x})| \geq 2$,
and let $I$ be a closed interval in $S^1$ of minimal length containing all the non-zero $x_{k}$. 
Then $0\in x_{1}\boxplus \cdots \boxplus x_{n}$ if and only if the length $\ell(I)$ of $I$ greater than or equal to $\pi$.
\end{lem}
\begin{proof}
Without loss of generality, suppose $x_1 \neq 0$.
This determines two closed half-circles, $C$ and $C'$, namely the two closed arcs joining $x_1$ and $-x_1$ both of length $\pi$.

It is not difficult to see that $\ell(I)=\pi$ if and only if 
\begin{enumerate}
\item at least one non-zero entry of $\mathbf{x}$ is $-x_1$ and 
\item every non-zero entry $x_k$ not equal to $x_1$ nor $-x_1$ is contained in the interior of exactly one of the half-circles $C$ or $C'$.
\end{enumerate}
If every non-zero entry of $\mathbf{x}$ is either $x_1$ or $-x_1$, then $I$ can be either $C$ or $C'$.
So assume that all of the non-zero entries of $\mathbf{x}$ are distinct and not antipodal.
Then all of the entries either lie in the interior of $C$ or of $C'$.
Without loss of generality, suppose $x_2$ and $x_3$ are the furthest point from $x_1$ in $C$ and $C'$, respectively.
Thus, $I$ is the union $x_1 \boxplus x_2$ and $x_1 \boxplus x_3$.
Then, $\ell(I) > \pi$ if and only if $-x_2 \in x_1 \boxplus x_3$ (or, equivalently,  $-x_3 \in x_1 \boxplus x_2$) if and only if $0 \in x_1 \boxplus x_2 \boxplus x_3$.
\end{proof}

The proof of this Lemma \ref{zero} implies the following.

\begin{lem}\label{pieces}
Let $(x_1,\ldots,x_n) \in \Phi^n-\{\mathbf{0}\}$ and $n \geq 3$.
If $0 \in x_1 \boxplus \cdots \boxplus x_n$, then there exists $j < k < \ell$ such that $0 \in x_j \boxplus x_k \boxplus x_\ell$.
\end{lem}

\subsection{Important subposets of $\mathbf{1}_n^\perp$}
We decompose $\mathbf{1}_n^\perp$ into the union of 
$$\{\mathbf{x} \in \mathbf{1}_n^\perp : x_n \in S^1\} \text{ and } \{\mathbf{x} \in \mathbf{1}_n^\perp : x_n=0\}.$$
Here we will address the former member of the union and assume $n \geq 3$. 
(The case when $n=2$ and the latter member of the union will be addressed in Section 4.)

By Lemma \ref{pieces}, we can write 
$$\{\mathbf{x} \in \mathbf{1}_n^\perp : x_n \in S^1\}=
\bigcup\limits_{j < k < \ell} \{\mathbf{x} \in \Phi^n-\{\mathbf{0}\}: 0 \in x_j \boxplus x_k \boxplus x_\ell
\text{ and } x_n \in S^1\}.$$
In fact, if $\mathbf{x} \in \Phi^n-\{\mathbf{0}\}$ such that $x_n \in S^1$ and $0 \in x_j \boxplus x_k \boxplus x_\ell$ for some $j < k < \ell <n$,
then either $0 \in x_j \boxplus x_\ell \boxplus x_n$ or $0 \in x_k \boxplus x_\ell \boxplus x_n$.
Thus we can write
$$\{\mathbf{x} \in \mathbf{1}_n^\perp : x_n \in S^1\}=
\bigcup\limits_{j < k < n} \{\mathbf{x} \in \Phi^n-\{\mathbf{0}\}: 0 \in x_j \boxplus x_k \boxplus x_n
\text{ and } x_n \in S^1\}.$$

Furthermore, suppose $\mathbf{x} < \mathbf{y}$ in $\{\mathbf{x} \in \mathbf{1}_n^\perp : x_n \in S^1\}$.
If $0 \in x_j \boxplus x_k \boxplus x_n$ for some $j < k <n$, then  $0 \in y_j \boxplus y_k \boxplus y_n$.
Thus every chain in $\{\mathbf{x} \in \mathbf{1}_n^\perp : x_n \in S^1\}$ is entirely contained in 
$\{\mathbf{x} \in \Phi^n-\{\mathbf{0}\}: 0 \in x_j \boxplus x_k \boxplus x_n\}$ for some $j < k < n$.
Hence
\begin{center}
$\Delta\{\mathbf{x} \in \mathbf{1}_n^\perp : x_n \in S^1\}=
\bigcup\limits_{j < k < n}\Delta\{\mathbf{x} \in \Phi^n-\{\mathbf{0}\}: 0 \in x_j \boxplus x_k \boxplus x_n\}$.
\end{center}

This observation together with Lemma \ref{homeo}, $\Delta\{\mathbf{x} \in \mathbf{1}_n^\perp : x_n \in S^1\}$ can be identified
with the union of 
$$\bigcup_{j < n} \{\mathbf{z} \in \mathbb{B}^n: (z_j,z_n) \in (S^1)^2 \text{ and } z_j=-z_n\}$$
and
$$\bigcup_{j<k<n} \{\mathbf{z} \in \mathbb{B}^n : (z_j,z_k,z_n) \in (S^1)^3 \text{ and } 0 \in z_j \boxplus z_k \boxplus z_n\}.$$
To understand how these pieces intersect, we will first analyze the case when $z_n=1$.

\subsection{Decomposition of $\Delta\{\mathbf{x} \in \mathbf{1}_n^\perp : x_n = 1\}$}
Given the contents of the previous subsection, we can identify $\Delta\{\mathbf{x} \in \mathbf{1}_n^\perp : x_n = 1\}$ with the union of
\begin{enumerate}
\item $\bigcup_{j < n} \{\mathbf{z} \in \mathbb{B}^n: (z_j,z_n)=(-1,1)\},$
\item $\bigcup_{j<k<n} \{\mathbf{z} \in \mathbb{B}^n : (z_j,z_k,z_n) \in U \times L \times \{1\} 
\text{ and } 0 \in z_j \boxplus z_k \boxplus 1\},$ and
\item $\bigcup_{j<k<n} \{\mathbf{z} \in \mathbb{B}^n : (z_j,z_k,z_n) \in L \times U \times \{1\} 
\text{ and } 0 \in z_j \boxplus z_k \boxplus 1\}.$
\end{enumerate}
Thus every member of the union in (1) is $\mathbb{B}^{n-2}$ which is homeomorphic to a $(2n-4)$-ball.
Since (2) can be written as $$\bigcup_{j<k<n} \{\mathbf{z} \in \mathbb{B}^n : (z_j,z_k,z_n) = (e^{i\pi t_j}, e^{i(\pi t_k +\pi)},1)
\text{ and } 0 \leq t_k  \leq t_j \leq 1\},$$ every member of this union is the product of a 2-ball and $\mathbb{B}^{n-3}$
which is homeomorphic to a $(2n-4)$-ball.

Similarly, (3) can be written as $$\bigcup_{j<k<n} \{\mathbf{z} \in \mathbb{B}^n : (z_j,z_k,z_n) \in (e^{i(\pi t_j +\pi)},e^{i\pi t_k},1)
\text{ and } 0 \leq t_j  \leq t_k \leq 1\},$$ a union of $(2n-4)$-balls.

Consider poset $P$ given by the Hasse diagram
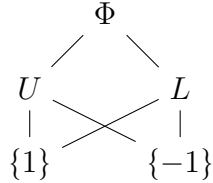
\begin{center}
\begin{tikzpicture}
\node (b) at (-1,1) {$\{1\}$};
\node (c) at (1,1) {$\{-1\}$};
\node (d) at (-1,2) {$U$};
\node (e) at (1,2) {$L$};
\node (f) at (0,3) {$\Phi$};

\draw (b) -- (d) -- (f); 
\draw (c) -- (e) -- (f); 
\draw (b) -- (e);
\draw (c) -- (d);
\end{tikzpicture}
\end{center}

Elements of $P^n$ will be treated as functions $[n] \to P$.
Define $\mathcal{P}_n \subset (P-\{0\})^n$ so that $X \in \mathcal{P}_n$ if and only if the following hold:
\begin{enumerate}
\item $X(\alpha)=\{1\}$ if and only if $\alpha=n$,
\item $X(\alpha) \in \{\{-1\}, U,L\}$ for some $\alpha \in [n-1]$, and 
\item $X(\alpha) \in \{U,L\}$ implies that there exists $\beta \in [n-1]-\{\alpha\}$ such that $X(\beta) \in \{U,L,-1\}$
but $X(\beta) \neq X(\alpha)$.
\end{enumerate}

\begin{example}
$(U,L,\Phi,1)$ and $(U,U,-1,1)$ are in $\mathcal{P}_4$ but not $(U,U,\Phi,1)$.
\end{example}

Let $X \in \mathcal{P}_n$ and define $B(X)$ to be the space of elements
$\mathbf{z} \in \mathbb{B}^n$ such that 
\begin{enumerate}
\item $z_\alpha \in \Delta(X(\alpha)) \text{ for each } \alpha \in [n]$, and 
\item $(z_\alpha,z_\beta) \in U \times L$ implies $(z_\alpha,z_\beta) = (e^{i\pi t_\alpha}, e^{i(\pi t_\beta +\pi)})$,
where $0 \leq t_\beta  \leq t_\alpha \leq 1$.
\end{enumerate}
Here we are identifying $\Delta(\Phi)$ with $\mathbb{B}$ and notice that the second condition implies that $0 \in z_\alpha \boxplus z_\beta \boxplus 1$.

An element $\mathbf{z} \in B(X)$ is contained in the \textit{interior} of $B(X)$, written as $\interior(B(X))$, if
\begin{enumerate}
\item $z_\alpha \in \interior(\Delta(X(\alpha))) \text{ for each } \alpha \in [n]$, and 
\item $(z_\alpha,z_\beta) \in U \times L$ implies $(z_\alpha,z_\beta) = \left(e^{i\pi t_\alpha}, e^{i(\pi t_\beta +\pi)}\right)$
such that $0 < t_\beta  < t_\alpha < 1$.
\end{enumerate}
Here $\interior(U)=\{e^{i\pi t}: 0<t<1\}$, $\interior(L)=\{e^{i(\pi t + \pi)}: 0<t<1\}$, $\interior(\{\text{point}\})=\{\text{point}\}$ and
$\interior(\mathbb{B})=\{re^{i2\pi}: r \in [0,1) \text{ and } t \in [0,1]\}$.

The \textit{boundary} of $B(X)$ is $\partial(B(X))=B(X)-\interior(B(X))$.

Define $\nu\colon \mathcal{P}_n \to [2,2n-4]_\mathbb{Z}$ by $\nu(X)=|X^{-1}(U)|+|X^{-1}(L)|+2|X^{-1}(\Phi)|$.
\begin{lem}\label{dimension}
$B(X)$ is a $\nu(X)$-ball.
\end{lem}
\begin{proof}
Clearly $B(X)$ is homeomorphic to
\begin{center}
$\{\mathbf{t} \in [0,1]^{X^{-1}(U) \cup X^{-1}(L)}: (x_\alpha, x_\beta) \in U \times L \text{ implies } t_\beta \leq t_\alpha\} 
\times \mathbb{B}^{X^{-1}(\Phi)}$.
\end{center}
Since the first member of the product is homeomorphic to a $(|X^{-1}(U)|+|X^{-1}(L)|)$-ball and 
the second is homeomorphic to a $(2|X^{-1}(\Phi)|)$-ball,
$B(X)$ is a $\nu(X)$-ball.
\end{proof}

\begin{notation}
Notice that $B(X)$ is a $(2n-4)$-ball if and only if $X$ is either of the form
\begin{center}
$X(\alpha)=\begin{cases}
-1 &\text{ if } \alpha=j,\\
1 &\text{ if } \alpha=n,\\
\Phi &\text{ otherwise}
\end{cases}$
\hspace{5pt}or\hspace{5pt}
$X(\alpha)=\begin{cases}
U &\text{ if } \alpha=j,\\
L &\text{ if } \alpha=k,\\
1 &\text{ if } \alpha=n,\\
\Phi &\text{ otherwise}
\end{cases}$
\end{center} 
for some integers $1 \leq j < k \leq n-1$. 
For convenience, we will denote $X^{(j,j)}$ for the former case and $X^{(j,k)}$ for the latter case.
\end{notation}

\begin{lem}\label{boundary}
Let $X$ and $Y$ be in $\mathcal{P}_n$.
If $X < Y$ in $\mathcal{P}_n$, then $B(X) \subset \partial(B(Y))$.
\end{lem}

\begin{proof}
Suppose $X < Y$ in $\mathcal{P}_n$ and $\mathbf{z} \in B(X)$.
Then, $X(\alpha) \leq Y(\alpha)$ for every $\alpha \in [n-1]$ which implies $z_\alpha \in \Delta(Y(\alpha))$.
If there exists an $(\alpha,\beta) \in X^{-1}(U) \times X^{-1}(L)$, then either
\begin{enumerate}
\item $X(\alpha) = Y(\alpha)$ and $X(\beta) = Y(\beta)$ 
which implies $0 \leq t_\beta  \leq t_\alpha \leq 1$ must be satisfied, or
\item $X(\alpha) < Y(\alpha)$ or $X(\beta) < Y(\beta)$ which implies  $0 \leq t_\beta  \leq t_\alpha \leq 1$ is not required.
\end{enumerate}
In either case, $\mathbf{z} \in B(Y)$. 
Furthermore, if $X(\alpha) < Y(\alpha)$, then we have the following cases:
\begin{enumerate}
\item $X(\alpha)=\{-1\}$ and $Y(\alpha) \in \{U,L,\Phi\}$ or
\item $X(\alpha) \in \{U,L\}$ and $Y(\alpha) = \Phi$.
\end{enumerate}
In either case, $z_\alpha \notin \interior(\Delta(Y(\alpha)))$.
Thus $\mathbf{z}$ cannot be contained in $\interior(B(Y))$.
\end{proof}

Clearly $\mathcal{P}_n$ is a lattice.

\begin{lem}\label{intersection}
Let $j \in [n-1]$ and $J \subseteq [n-1]$ such that $|J|>1$. 
Then 
$$\bigcap\limits_{k \in J}B\left(X^{(j,k)}\right)=B\left(\bigwedge\limits_{k \in J}X^{(j,k)}\right) \text{ and }
\bigcap\limits_{k \in J}B\left(X^{(k,j)}\right)=B\left(\bigwedge\limits_{k \in J}X^{(k,j)}\right).$$
\end{lem}
The proof follows from Lemma \ref{boundary} and definition. 

\begin{theorem}\label{slice}
$\Delta\{\mathbf{x} \in \mathbf{1}_n^\perp : x_n = 1\}$ is a $(2n-4)$-ball. 
\end{theorem}
\begin{proof}
For every $j \in [n-1]$, define $\mathcal{B}_j \coloneqq \bigcup\limits_{k \in [n-1]} B\left(X^{(j,k)}\right)$.
Therefore
\begin{center}
$\Delta\{\mathbf{x} \in \mathbf{1}_n^\perp : x_n = 1\}=\bigcup\limits_{j \in [n-1]} \mathcal{B}_j$.
\end{center}
To apply Theorem \ref{gluing} on $\{\mathcal{B}_j: j \in [n-1]\}$, fix $j \in [n-1]$.

\begin{cl}
$\mathcal{B}_j$ is a  $(2n-4)$-ball.
\end{cl}
\begin{proof}
Consider $\left\{X^{(j,k)}: k \in [n-1]\right\}$.
Since $\nu\left(X^{(j,k)}\right)=2n-4$ for every $k \in [n-1]$, $\mathcal{B}_j$ is a union of $(2n-4)$-balls by Lemma \ref{dimension}.

Now fix $J \subseteq [n-1]$, where $|J|>1$. 
Since $\bigcap\limits_{k \in J}B\left(X^{(j,k)}\right)=B\left(\bigwedge\limits_{k \in J}X^{(j,k)}\right)$ by Lemma \ref{intersection}, 
we have $\nu\left(\bigwedge\limits_{k \in J}X^{(j,k)}\right)=2n-3-|J|$.
Thus  $\bigcap\limits_{k \in J}B\left(X^{(j,k)}\right)$ is $(2n-3-|J|)$-ball by Lemma \ref{dimension}.

Lastly, fix $r \in J$. 
Since $\bigwedge\limits_{k \in J}X^{(j,k)} < \bigwedge\limits_{k \in J-\{r\}}X^{(j,k)}$ and $B\left(\bigwedge\limits_{k \in J-\{r\}}X^{(j,k)}\right)=\bigcap\limits_{k \in J-\{r\}}B\left(X^{(j,k)}\right)$ by Lemma \ref{intersection}, we have 
\begin{center}
$\bigcap\limits_{k \in J}B\left(X^{(j,k)}\right) \subseteq \partial\left(\bigcap\limits_{k \in J-\{r\}}B\left(X^{(j,k)}\right)\right)$
\end{center}
by Lemma \ref{boundary}.
Thus, by Theorem \ref{gluing}, $\mathcal{B}_j$ is a $(2n-4)$-ball.
\end{proof}

To show the first condition of Theorem \ref{gluing}, fix $J \subseteq [n-1]$ such that $|J|>1$.
\begin{cl}
$\bigcap\limits_{j \in J} \mathcal{B}_j= \bigcup\limits_{k \in [n-1]} B\left(\bigwedge\limits_{j \in J} X^{(j,k)}\right)$ and is a $(2n-3-|J|)$-ball.
\end{cl}
\begin{proof}
To prove the first part of the claim, we will induct on the cardinality of $J$.

When $|J|=2$, without loss of generality, assume $J=\{1,2\}$.
By definition
\begin{center}
\begin{align*}
\mathcal{B}_1 \cap \mathcal{B}_2
&=\left(\bigcup_{k \in [n-1]}B\left(X^{(1,k)}\right)\right) \cap \left(\bigcup_{l \in [n-1]}B\left(X^{(2,\ell)}\right)\right)\\
&=\bigcup_{k \in [n-1]}\bigcup_{\ell \in [n-1]} B\left(X^{(1,k)}\right) \cap B\left(X^{(2,\ell)}\right)
\end{align*}
\end{center}

If $k=\ell$, then $B\left(X^{(1,k)}\right) \cap B\left(X^{(2,\ell)}\right)=B\left(X^{(1,k)} \wedge X^{(2,k)}\right)$ by Lemma \ref{intersection}.

If $k\neq\ell$, then we must consider the cases when (1)  $k\in[2]$, (2) $\ell \in [2]$, or (3) neither $k$ nor $\ell$ are in $[2]$.
In cases (1) and (2), it is easy to show that $B\left(X^{(1,k)}\right) \cap B\left(X^{(2,\ell)}\right)$ the union of 
$B\left(X^{(1,2)} \wedge X^{(2,2)}\right)$ and $B\left(X^{(1,1)} \wedge X^{(2,1)}\right)$.
In case (3), let $\mathbf{z} \in B\left(X^{(1,k)}\right) \cap B\left(X^{(2,\ell)}\right)$. 
Then $(z_1,z_k,z_2,z_\ell)=(e^{i\pi t_1},e^{i(\pi t_k +\pi)},e^{i\pi t_2},e^{i(\pi t_\ell + \pi)})$, 
where $0 \leq t_k \leq t_1 \leq 1$ and $0 \leq t_\ell \leq t_2 \leq 1$.
If $t_k \leq t_\ell$, then $\mathbf{z} \in B\left(X^{(1,\ell)} \wedge X^{(2,\ell)}\right)$.
If $t_\ell \leq t_k$, then $\mathbf{z} \in B\left(X^{(1,k)} \wedge X^{(1,k)}\right)$.
Thus 
\begin{center}
$\mathcal{B}_1 \cap \mathcal{B}_2=\bigcup\limits_{k \in [n-1]} B\left(X^{(1,k)} \wedge X^{(2,k)}\right)$.
\end{center}

For the induction step, suppose that the claim is true for $J=[m]$, where $3 \leq m \leq n-1$, i.e., 
\begin{center}
$\bigcap\limits_{j \in [m]}\mathcal{B}_j= \bigcup\limits_{k \in [n-1]} B\left(\bigwedge\limits_{j \in [m]} X^{(j,k)}\right)$.
\end{center}
Thus
\begin{center}
\begin{align*}
\bigcap\limits_{j \in [m+1]}\mathcal{B}_j &=\mathcal{B}_m \cap \left(\bigcap\limits_{j \in [m]}\mathcal{B}_j\right)\\
&=\left(\bigcup_{\ell \in [n-1]} B\left(X^{(m+1,\ell)}\right)\right) \cap \left(\bigcup\limits_{k \in [n-1]} B\left(\bigwedge\limits_{j \in [m]} X^{(j,k)}\right)\right)\\
&=\bigcup\limits_{\ell \in [n-1]}\bigcup\limits_{k \in [n-1]} \left(B\left(X^{(m+1,\ell)}\right) \cap B\left(\bigwedge\limits_{j \in [m]} X^{(j,k)}\right)\right).
\end{align*}
\end{center}

If $k=\ell$, then 
$B\left(X^{(m+1,\ell)}\right) \cap B\left(\bigwedge\limits_{j \in [m]} X^{(j,k)}\right)=B\left(\bigwedge\limits_{j \in [m+1]} X^{(j,k)}\right)$
by Lemma \ref{intersection}.

If $k\neq\ell$, then we must consider the cases when 
\begin{enumerate}
\item $k\in[m]$ which implies $B\left(X^{(m+1,\ell)}\right) \cap B\left(\bigwedge\limits_{j \in [m]} X^{(j,k)}\right)$ is contained in $B\left(\bigwedge\limits_{j \in [m+1]} X^{(j,k)}\right)$, 
\item $\ell \in [m]$ which implies $B\left(X^{(m+1,\ell)}\right) \cap B\left(\bigwedge\limits_{j \in [m]} X^{(j,k)}\right)$ is contained in $B\left(\bigwedge\limits_{j \in [m+1]} X^{(j,\ell)}\right)$, 
\item neither $k$ nor $\ell$ are in $[m]$.
\end{enumerate}
In case (3), let $\mathbf{z} \in B(X^{(m+1,\ell)}) \cap B\left(\bigwedge\limits_{j \in [m]} X^{(j,k)}\right)$
which implies $(z_\ell,z_k)=\left(e^{i(\pi t_\ell + \pi)},e^{i(\pi t_\ell + \pi)}\right)$, where $t_\ell,t_k \in [0,1]$. 
Then, either $t_\ell \leq t_k$ which implies $\mathbf{z} \in B\left(\bigwedge\limits_{j \in [m+1]} X^{(j,k)}\right)$,
or $t_k \leq t_\ell$ which implies $\mathbf{z} \in B\left(\bigwedge\limits_{j \in [m+1]} X^{(j,\ell)}\right)$.
Thus
$\bigcap\limits_{j \in J}\mathcal{B}_j= \bigcup\limits_{k \in [n-1]} B\left(\bigwedge\limits_{j \in J} X^{(j,k)}\right)$.

To prove that the second part of the claim, we are going to apply Theorem \ref{gluing} on  the collection $\left\{B\left(\bigwedge\limits_{j \in J} X^{(j,k)}\right) : k \in [n-1]\right\}$.
Clearly, $\nu\left(\bigwedge\limits_{j \in J} X^{(j,k)}\right)=2n-3-|J|$ for each $k\in [n-1]$.
Thus $B\left(\bigwedge\limits_{j \in J} X^{(j,k)}\right)$ is a $(2n-3-|J|)$-ball by Lemma \ref{dimension}.

To show that the first condition of Theorem \ref{gluing} holds, fix $I \subseteq [n-1]$, where $|I|>1$.
It is not difficult to see that 
\begin{center}
$\bigcap\limits_{k \in I}B\left(\bigwedge\limits_{j \in J} X^{(j,k)}\right)=B\left(\bigwedge\limits_{(j,k) \in J \times I} X^{(j,k)}\right)$
\end{center}
and $\nu\left(\bigwedge\limits_{(j,k) \in J \times I} X^{(j,k)}\right)=2n-2-|J|-|I|$.
Hence, $\bigcap\limits_{k \in I}B\left(\bigwedge\limits_{j \in J} X^{(j,k)}\right)$ is a $(2n-2-|J|-|I|)$-ball by Lemma \ref{dimension}.

To show that the second condition holds of Theorem \ref{gluing}, fix $r \in I$.
We must show that
\begin{center}
$\bigcap\limits_{k \in I-\{r\}} B\left(\bigwedge\limits_{j \in J} X^{(j,k)}\right) \subseteq \partial\left(\bigcap\limits_{k \in I} B\left(\bigwedge\limits_{j \in J} X^{(j,k)}\right)\right)$
\end{center}
which follows immediately from Lemma \ref{intersection} and Lemma \ref{boundary}.
Thus, by Theorem \ref{gluing}, 

 \noindent $\bigcup\limits_{k \in [n-1]} B\left(\bigwedge\limits_{j \in J} X^{(j,k)}\right)$ is a $(2n-3-|J|)$-ball.
\end{proof}

To show the second condition of Theorem \ref{gluing} holds for $\{\mathcal{B}_j: j \in [n-1]\}$, fix $r \in J$.
\begin{cl}
$\bigcap\limits_{j \in J} \mathcal{B}_j \subseteq \partial\left(\bigcap\limits_{j \in J-\{r\}} \mathcal{B}_j\right)$.
\end{cl}
\begin{proof}
We must show that
\begin{center}
$\bigcup\limits_{k \in [n-1]} B\left(\bigwedge\limits_{j \in J} X^{(j,k)}\right)=\bigcap\limits_{j \in J} \mathcal{B}_j \subseteq \partial\left(\bigcap\limits_{j \in J-\{r\}} \mathcal{B}_j\right)=\partial\left(\bigcup\limits_{k \in [n-1]} B\left(\bigwedge\limits_{j \in J} X^{(j,k)}\right)\right)$.
\end{center}
Fix $k_0 \in [n-1]$.
It suffices to show that $B\left(\bigwedge\limits_{j \in J} X^{(j,k_0)}\right) \subseteq \partial\left(\bigcup\limits_{k \in [n-1]} B\left(\bigwedge\limits_{j \in J} X^{(j,k)}\right)\right)$
which follows from Lemma \ref{intersection} and Lemma \ref{boundary}.
Thus $\bigcap\limits_{j \in J} \mathcal{B}_j \subseteq \partial\left(\bigcap\limits_{j \in J-\{r\}} \mathcal{B}_j\right)$.
\end{proof}
Hence, by Theorem \ref{gluing}, $\Delta\{\mathbf{x} \in \mathbf{1}_n^\perp: x_n=1\}$ is a $(2n-4)$-ball.
\end{proof}

\section{Proof of Theorem \ref{main}}

We will induct on $n \geq 2$.

When $n=2$, we have $(1,1)^\perp=\{(0,0)\} \cup \{(x,-x): x \in S^1\}$.
Since  $(1,1)^\perp-\{(0,0)\}$ contains no chains, $\Delta((1,1)^\perp-\{(0,0)\})=(1,1)^\perp-\{(0,0)\}$
which is clearly homeomorphic to $S^1$.

For $n>2$, assume that the Theorem \ref{main} holds for $n-1$, that is, we have a homeomorphism $h \colon \Delta\left(\mathbf{1}_{n-1}^\perp-\{\mathbf{0}\}\right) \to S^{2n-5}.$
As previously discussed, we can write 
\begin{center}
$\mathbf{1}_n^\perp=\{\mathbf{x} \in \mathbf{1}_n^\perp: x_n \in S^1\} \cup \{\mathbf{x} \in\mathbf{1}_n^\perp : x_n=0\}.$
\end{center}
By Theorem \ref{slice}, $\Delta(\{\mathbf{x} \in \mathbf{1}_n^\perp : x_n = 1\})$ is a $(2n-4)$-ball.
\begin{cl}
$\partial\left(\Delta(\{\mathbf{x} \in \mathbf{1}_n^\perp : x_n = 1\})\right)=
\Delta\{(\mathbf{y},1): \mathbf{y} \in \mathbf{1}_{n-1}^\perp-\{\mathbf{0}\}\}$.
\end{cl}
\begin{proof}
First we show that
$\Delta\{(\mathbf{y},1): \mathbf{y} \in \mathbf{1}_{n-1}^\perp-\{\mathbf{0}\}\} \subseteq 
\partial\left(\Delta\left\{\mathbf{x} \in \mathbf{1}_n^\perp: x_n=1\right\}\right)$.
Let $\mathbf{z} \in \Delta\{(\mathbf{y},1): \mathbf{y} \in \mathbf{1}_{n-1}^\perp-\{\mathbf{0}\}\} 
\subseteq \Delta\left\{\mathbf{x} \in \mathbf{1}_n^\perp: x_n=1\right\}$.
Since $\mathbf{z} \in \Delta\left\{\mathbf{x} \in \mathbf{1}_n^\perp: x_n=1\right\}$, 
we have that $\mathbf{z} \in B\left(\bigwedge\limits_{j \in J} X^{(j,k)}\right)$ for some $J \subseteq [n-1]$ and $k \in [n-1]$.
To show that $\mathbf{z} \in \partial\left(\Delta\left\{\mathbf{x} \in \mathbf{1}_n^\perp: x_n=1\right\}\right)$,
it suffices to show that  $\mathbf{z} \notin \interior\left(B\left(\bigwedge\limits_{j \in J} X^{(j,k)}\right)\right)$.
By way of contradiction, suppose that $\mathbf{z} \in \interior\left(B\left(\bigwedge\limits_{j \in J} X^{(j,k)}\right)\right)$, then
it has to satisfy one of two lists of conditions depending on $k$.

If $k \in J$, then
\begin{enumerate}
\item $z_k=-1$,
\item $z_j \in \interior U$ for every $j  \in J-\{k\}$, and
\item $z_\ell \in \interior \mathbb{B}$ for $\ell \in [n-1]-J$.
\end{enumerate}

If $k \in [n-1]-J$, then
\begin{enumerate}
\item $z_k=e^{i(\pi t_k + \pi)}$ and $z_j=e^{i \pi t_j}$, where $0< t_j < t_k$  for every $j  \in J$, and 
\item $z_\ell \in \interior\mathbb{B}$ for $\ell \in [n-1]-J$.
\end{enumerate}

In either case, $0 \notin z_k \bigboxplus_{j \in J} z_j$, a contradiction.
Thus $\mathbf{z} \in \interior\left(B\left(\bigwedge\limits_{j \in J} X^{(j,k)}\right)\right)$ and, hence, 
\begin{center}
$\Delta\{(\mathbf{y},1): \mathbf{y} \in \mathbf{1}_{n-1}^\perp-\{\mathbf{0}\}\} \subseteq \partial\left(\Delta\left\{\mathbf{x} \in \mathbf{1}_n^\perp: x_n=1\right\}\right)$.
\end{center}

To obtain equality, it suffices to show that 
$\Delta\{(\mathbf{y},1): \mathbf{y} \in \mathbf{1}_{n-1}^\perp-\{\mathbf{0}\}\}$ 
is homeomorphic to $S^{2n-5}$.
Consider the map
\begin{center}
$\pi_{[n-1]}\colon \Phi^n \to \Phi^{n-1},\ \mathbf{x} \mapsto (x_1,x_2,\ldots,x_{n-1})$.
\end{center}
We can restrict the domain and the codomain of $\pi_{[n-1]}$ to obtain
\begin{center}
$\phi\colon \{(\mathbf{y},1):\mathbf{y} \in \mathbf{1}_{n-1}^\perp-\{\mathbf{0}\}\} \to \mathbf{1}_{n-1}^\perp-\{\mathbf{0}\}$.
\end{center}
This map is an isomorphism of topological posets which induces a homeomorphism $\Delta(\phi)$ between their topological order complexes.
By the induction hypothesis of Theorem \ref{main}, we have that
\begin{center}
$h \circ \Delta(\phi)\colon
\Delta\{(\mathbf{y},1):  \mathbf{1}_{n-1}^\perp-\{\mathbf{0}\}\} \to S^{2n-5}$
\end{center}
is a homeomorphism.
Thus,
\begin{center}
$\Delta\{(\mathbf{y},1):  \mathbf{1}_{n-1}^\perp-\{\mathbf{0}\}\}=\partial\left(\Delta\left\{\mathbf{x} \in \mathbf{1}_n^\perp: x_n=1\right\}\right)$.
\end{center}
\end{proof}

We can extend $h \circ \Delta(\varphi)$ to a homeomorphism 
$H \colon \Delta\left\{\mathbf{x} \in \mathbf{1}_n^\perp: x_n=1\right\} \to \mathbb{B}^{n-2}$.
Since $(-1,\dots,-1,1) \in \interior\left(\Delta\left\{\mathbf{x} \in \mathbf{1}_n^\perp: x_n=1\right\}\right)$, we can define the unique path 
\begin{center}
$p_\mathbf{z}\colon[0,1] \to \mathbb{B}^{n-2}$\\
$t \mapsto (1-t)H(-1,\dots,-1,1)+tH(\mathbf{z})$
\end{center}
for every $\mathbf{z} \in \partial\left(\Delta\left\{\mathbf{x} \in \mathbf{1}_n^\perp-\{\mathbf{0}\}: x_n=1\right\}\right)$.
Thus, we have a unique path from $(-1,\dots,-1,1)$ to $\mathbf{z}$ via $H^{-1} \circ p_\mathbf{z}$.
Furthermore, there is a unique path 
\begin{center}
$p_\mathbf{z}'\colon[0,1] \to \Delta(\mathbf{1}_n^\perp-\{0\}),\ t \mapsto (z_1,z_2,\dots,(1-t)z_n)$.
\end{center}
from every $\mathbf{z} \in \partial\left(\Delta\left\{\mathbf{x} \in \mathbf{1}_n^\perp-\{\mathbf{0}\}: x_n=1\right\}\right)$
to $\mathbf{z}' \in \Delta\{\mathbf{x} \in\mathbf{1}_n^\perp : x_n=0\}$, where $\mathbf{z}'=\mathbf{z}$ except $z_n'=0$.
Therefore,
\begin{center}
$[0,1] \to \Delta(\mathbf{1}_n^\perp-\{0\}),\  \mapsto
\begin{cases}
H^{-1} \circ p_\mathbf{z}(2t) &\text{if } t \in [0,1/2],\\
p_\mathbf{z}'(2t-1) &\text{if } t \in [1/2,1].
\end{cases}
$
\end{center}
is a unique path from $(-1,-1,\dots,-1,1)$ to every $\mathbf{z}' \in \Delta\{\mathbf{x} \in\mathbf{1}_n^\perp : x_n=0\}$.
Since, for every $x \in S^1$, $g_x|\colon \Delta(\mathbf{1}_n^\perp-\{\mathbf{0}\}) \to \Delta(\mathbf{1}_n^\perp-\{\mathbf{0}\})$
is a homeomorphism, we have a unique path 
from $(-x,-x,\dots,-x,x) \in g_x (\Delta\{\mathbf{x} \in \mathbf{1}_n^\perp : x_n = 1\})$ 
to every point in $\Delta\{\mathbf{x} \in\mathbf{1}_n^\perp : x_n=0\}$.
In other words, $\Delta(\mathbf{1}_n^\perp-\{\mathbf{0}\})$ is the join of the circle
$\{(-x,-x,\dots,-x,x): x \in S^1\}$ and the $(2n-5)$-sphere $\Delta\{\mathbf{x} \in\mathbf{1}_n^\perp : x_n=0\}$.
Thus  $\Delta(\mathbf{1}_n^\perp-\{\mathbf{0}\})$ is homeomorphic to $S^1 \ast S^{2n-5} \cong S^{2n-3}.$

\bibliographystyle{amsalpha}
\bibliography{sources}{} 

\end{document}